\documentclass{amsart}
\usepackage[latin1]{inputenc}
\newtheorem{thm}{Theorem}
\newtheorem{cor}[thm]{Corollary}
\newtheorem{lem}[thm]{Lemma}
\newtheorem{prop}[thm]{Proposition}

\theoremstyle{definition}
\newtheorem{defn}[thm]{Definition}

 \numberwithin{equation}{section}
\newcommand{\To}{\longrightarrow}

\begin{document}

\thanks{This research was partially supported by the grant
BFM2002-01719 of MCyT (Spain) and a FPU grant of MEC (Spain).}
 %\keywords{weakly compact, index of $mathcal{K}$-determination}
%\date{\today}
\title[]{The number of weakly compact sets which generate a Banach space}%
\author{Antonio Avilés}

\address{Departamento de Matemáticas\\ Universidad de Murcia\\ 30100 Espinardo (Murcia)\\ Spain }%
\email{avileslo@um.es}

\begin{abstract}
We consider the cardinal invariant $CG(X)$ of the minimal number of weakly
compact subsets which generate a Banach space $X$. We study the behavior
of this index when passing to subspaces, its relation with the Lindel\"of
number in the
weak topology and other related questions.\\
\end{abstract}

\maketitle

A Banach space is weakly compactly generated if there is a weakly
compact subset which is linearly dense and weakly Lindelöf if it
is a Lindelöf space in its weak topology. It was asked by Corson
\cite{Corsonquestion} which was the relation between these two
concepts. The answer was that every weakly compactly generated
space is weakly Lindelöf but the converse is not true, and in
order to clarify what was in the middle the class of weakly
$\mathcal{K}$-analytic was introduced by Talagrand~\cite{WKA}, who
was together with Pol~\cite{PolCXLindelof} the first to solve this
problem. Here we shall analyze the question of Corson from a more
general point of view: What is the relation between the number of
weak compacta which are necessary to generate a Banach space and
the Lindelöf number of the space in the weak topology? Again, an
intermediate class analagous to that introduced by Talagrand plays
a clarifying role in the theory. Thus, our starting point is the
following (cf. Sections \ref{sectcard} and \ref{sectusco} for
notation):

\begin{defn}\label{definicionBanach}
Let $X$ be a topological space.
\begin{enumerate}
\item The index of compact generation of $X$, $CG(X)$, is defined
as the least infinite cardinal $\kappa$ such that there exists a family
$\{K_\lambda : \lambda<\kappa\}$ of compact subsets of $X$ whose union is
a dense subset of $X$.

\item The index of $\mathcal{K}$-analyticity of $X$, $\ell K(X)$,
is the least infinite cardinal $\kappa$ for which there exists a complete
metric space $M$ of weight $\kappa$ and an usco $M\To 2^{X}$.

\item The Lindelöf number of $X$, $\ell(X)$, is the least
infinite cardinal $\kappa$ such that any cover of $X$ by open sets has a
subcover with at most $\kappa$ many sets.
\end{enumerate}
\end{defn}

If $X$ is a Banach space, all the indices will refer always to the weak
topology of $X$. In this way the classes of weakly compactly generated,
weakly $\mathcal{K}$-analytic and weakly Lindelöf Banach spaces equal the
classes of spaces $X$ such that $CG(X)=\omega$, $\ell K(X)=\omega$ and
$\ell(X) =\omega$ respectively. Similar indices to $\ell K(X)$ can be
defined if instead of complete metric spaces of a given weight we use
other classes of topological spaces. These kind of indices have been
studied in \cite{Maria}, cf. \cite{CasMunOri}, such as the index of
$\mathcal{K}$-determinacy $\ell\Sigma(X)$ (taking in (2) arbitrary metric
spaces of weight $\kappa$ instead of complete metric spaces) and the
Nagami index $Nag(X)$ (taking arbitrary completely regular topological
spaces of weight $\kappa$). For any Banach space $X$ we have that

$$\ell(X) \leq Nag(X) \leq \ell\Sigma(X)\leq \ell K(X) \leq CG(X)$$

The first inequality can be found in \cite{CasMunOri} and follows from the
fact that the Lindelöf number cannot increase by an usco image. The others
are self-evident excepts perhaps the last one, cf. Section
\ref{sectrelacionCGLK} below. This gives a first answer to our original
question: the number of weak compacta which are necessary to generate a
Banach space is not lower than the Lindelöf number in the weak topology,
shortly $\ell(X)\leq CG(X)$. On the other hand, we shall show that $CG(X)$
may be arbitrarily larger than $\ell(X)$, namely in Section \ref{largeCG}
we prove:

\begin{thm}\label{bigCG}
Let $\kappa$ be any cardinal. There exists a weakly Lindelöf determined
Banach space $X$ such that $CG(X)>\kappa$.
\end{thm}

Weakly Lindelöf determined spaces constitute a special class of weakly
Lindelöf Banach spaces \cite{ArgMerWLD}. The relation of $CG(X)$ with the
other indices is quite different because the cardinality of a completely
regular space of weight $\kappa$ is at most $2^\kappa$, so $CG(X)\leq
2^{Nag(X)}$, and by the same reason
$CG(X)\leq \ell\Sigma(X)^\omega$.\\

Another classical question is the fact, first shown by
Rosenthal~\cite{SubspacesWCG}, that there are subspaces of weakly
compactly generated spaces which are not weakly compactly generated. It
has been observed in \cite{cardinalb} that such spaces must have large
enough density character, namely greater than or equal to cardinal
$\mathbf{b}$. In Section \ref{sectsubespacios}, we address the natural
question now: What is the relation between $CG(X)$ and $CG(Y)$ for $Y$ a
subspace of $X$? and what about the density character? The answer we give
is the following:

\begin{thm}\label{thCGsubespacios}
Let $\kappa,\tau,\delta$ be infinite cardinals. The following are
equivalent:
\begin{enumerate}
\item $\tau\leq\mathbf{d}(\kappa)$ and
$\delta\geq\mathbf{b}_\kappa(\tau)$.
\item There exists a Banach
space $X$ and a subspace $Y$ of $X$ such that $CG(X) = \kappa$, $CG(Y) =
\tau$ and $dens(Y)=\delta$.
\end{enumerate}
\end{thm}

The cardinal numbers $\mathbf{d}(\kappa)$ and $\mathbf{b}_\kappa(\tau)$
are defined in Section \ref{sectcard} in terms of the topology of the
space $\kappa^\omega$. For cardinals $\kappa,\tau\geq 2^\omega$, it
happens that $\mathbf{d}(\kappa)=\kappa^\omega$ and
$\mathbf{b}_\kappa(\tau)=\tau$ but for cardinals below the continuum the
behavior of this functions is more complicated and depends heavily on the
axiomatic settlement. The fact that (1) implies (2) in Theorem
\ref{thCGsubespacios} is obtained by modifying an example of Argyros
\cite[Section 1.6]{FabianWA} while the converse is based on the use of the
index $\mathcal{K}$-analyticity. It is indeed established in Section
\ref{sectrelacionCGLK} a very similar result to Theorem
\ref{thCGsubespacios} concerning the relation of the indices $\ell K(X)$
and $CG(X)$:

\begin{thm}\label{relacionCGWK}
Let $\kappa,\tau,\delta$ be infinite cardinals. The following are
equivalent:
\begin{enumerate}
\item $\kappa\leq\tau\leq\mathbf{d}(\kappa)$ and
$\delta\geq\mathbf{b}_\kappa(\tau)$.
\item There exists a Banach
space $X$ such that $\ell K(X) = \kappa$, $CG(X) = \tau$ and
$dens(X)=\delta$.
\end{enumerate}
\end{thm}

Here, the fact that (1) implies (2) is obtained by modifying the
construction of Talagrand \cite{WKA} of a weakly $\mathcal{K}$-analytic
space which is not weakly compactly generated. It remains unclear to us
what are the precise relations between the indices $\ell K(X)$,
$\ell\Sigma(X)$ and $Nag(X)$ for Banach spaces. We make some remarks about
that in Section~\ref{otherindices}.\\

We want to express our gratitude to Bernardo Cascales and José Orihuela
for their help and support during this work, as well as to David Fremlin
for valuable comments and suggestions.

\section{Cardinal numbers and metric spaces}\label{sectcard}

In this section we fix the notation about cardinal arithmetic and metric
spaces, and we shall define and discuss some cardinal numbers which will
be used in Sections \ref{sectrelacionCGLK} and \ref{sectsubespacios}. A
cardinal number $\kappa$ is identified with the set of all ordinals less
than $\kappa$, and in particular $\kappa$ is a set of cardinality
$\kappa$, and is considered also as a topological space endowed with the
discrete topology. By $\kappa^\omega$ we denote the set of all sequences
of elements of $\kappa$ endowed with the product topology (with respect to
the discrete topology on each factor), as well as the cardinality of this
set. Finally, $2^A$ stands for the family of all subsets of $A$, and when
$A=\kappa$ is a cardinal, $2^\kappa$ also denotes the cardinality of this
set.

\begin{defn}
Let $\kappa$ be a cardinal number:
\begin{enumerate}
\item The cardinal $\mathbf{d}(\kappa)$ is defined as the least
cardinal $\lambda$ such that $\kappa^\omega$ is the union of $\lambda$
many compact subsets.
\item Let $\tau$ be a cardinal such that
$\tau\leq\mathbf{d}(\kappa)$. The cardinal $\mathbf{b}_\kappa(\tau)$ is
the least cardinal $\lambda$ for which there exists a set $A$ of
cardinality $\lambda$ such that $A$ is not contained in any union of less
than $\tau$ many compact subsets of $\kappa^\omega$.
\end{enumerate}
\end{defn}

Notice that always $\kappa\leq\mathbf{d}(\kappa)\leq\kappa^\omega$
and $\tau\leq \mathbf{b}_\kappa(\tau)\leq\kappa^\omega$, and that
$\mathbf{b}_\kappa(\tau)=\tau$ if $\tau\leq\kappa$ (we can
consider $A$ a closed and discrete subset of $\kappa^\omega$ of
cardinality $\tau$). Every compact metrizable space is either
countable or has cardinality $2^\omega$, hence
$\mathbf{d}(\kappa)=\kappa^\omega$ whenever $\kappa>2^\omega$ and
$\mathbf{b}_\kappa(\tau)=\tau$ whenever
$2^\omega<\tau\leq\mathbf{d}(\kappa)$. On the other hand, if
$cof(\kappa)>\omega$ then $\kappa^\omega =
\bigcup_{\alpha<\kappa}\alpha^\omega$, and this fact implies that
$\mathbf{d}(\kappa)=\sum_{\alpha<\kappa}\mathbf{d}(|\alpha|)$. The
difficult case in computing $\mathbf{d}(\kappa)$ is when $\kappa$
is a cardinal of cofinality $\omega$ less than $2^\omega$. For
example, when $\kappa=\omega$ we refer to \cite{vanDouwen} for
information about cardinal $\mathbf{d} = \mathbf{d}(\omega)$. We
illustrate also the situation for $\kappa=\omega_\omega$, for
which we need the following observation, pointed out to us by
David Fremlin:

\begin{prop}\label{Fremlinsobs}
For an infinite cardinal $\kappa$, $\mathbf{d}(\kappa)=\max[\mathbf{d},
\mathrm{cf}([\kappa]^{\leq\omega})]$, where
$\mathrm{cf}([\kappa]^{\leq\omega})$ is the least cardinality of a cofinal
family $A$ of countable subsets of $\kappa$, that is, a family such that
every countable subset of $\kappa$ is a subset of some member of $A$.
\end{prop}

We mention that Shelah has obtained that
$\mathrm{cf}([\omega_\omega]^{\leq\omega})<\omega_{\omega_4}$, cf.
\cite{BurMagdeShelah}. The proof of Proposition \ref{Fremlinsobs}
is not difficult: if $B$ is a family of compact sets covering
$\kappa^\omega$, then the family $A = \{\{x_n :
x=(x_i)_{i<\omega}\in K, n<\omega\} : K\in B\}$ is a cofinal
family of countable subsets of $\kappa$, and conversely if $A$ is
a cofinal family of countable subsets of $\kappa$ and for every
$s\in A$, $C_s$ is a family of $\mathbf{d}$ many compact sets
covering $s^\omega$, then $B = \bigcup_{s\in A}C_s$ is a family of
compact sets
covering $\kappa^\omega$.\\

About cardinals $\mathbf{b}_\kappa(\tau)$ we know very little more except
for $\mathbf{b} = \mathbf{b}_\omega(\omega_1)$ \cite{vanDouwen} and the
fact that in some cases we can establish a relation with the known
$\mathbf{b}$, for example $\mathbf{b}\leq \mathbf{b}_\omega(\omega_2) =
\mathbf{b}_{\omega_1}(\omega_2)$ provided $\omega_1<\mathbf{d}$.

\section{The $\mathcal{K}$-analyticity index}\label{sectusco}

If $\Sigma$ and $Y$ are topological spaces, we will say that a map
$\phi:\Sigma\To 2^Y$ is an usco if the three following conditions are
satisfied:

\begin{enumerate}
\item $\phi(\sigma)$ is a compact subset of $Y$, for every $\sigma\in\Sigma$.
\item For every open subset $U$ of $Y$, the set $\{\sigma\in\Sigma :
\phi(\sigma)\subseteq U\}$ is an open subset of $\Sigma$.
\item $\bigcup_{\sigma\in\Sigma}\phi(\sigma) = Y$.
\end{enumerate}

In this situation, if $A$ is a subset of $\Sigma$, we denote
$\phi(A)=\bigcup\{\phi(\sigma) : \sigma\in A\}$.\\

We recall that any complete metric space $X$ of weight $\kappa$ is
a continuous image of a closed subset $M$ of $\kappa^\omega$: One
considers a base $\{O_\lambda : \lambda<\kappa\}$ of $X$ and $M =
\{x\in\kappa^\omega : diam(O_{x_n})<\frac{1}{n} ,
\overline{O_{x_{n+1}}}\subset O_{x_n}\}$. In addition, for $M$ a
closed subset of $\kappa^\omega$ there is retract
$p:\kappa^\omega\To M$ \cite[Proposition 2.8]{Kechrisbook}. In
particular, $\ell K(Y)\leq\kappa$ if and only if there is an usco
$\kappa^\omega\To 2^Y$. The following two properties that we shall
use can be found in \cite{CasMunOri} proven for the index
$\ell\Sigma(X)$ but the proof for $\ell K(X)$ is completely
analogous (as usual, $C_p(K)$ stands for the space of continuous
functions over $K$ endowed with the pointwise convergent
topology):

\begin{prop}\label{propiedadesMaria} Let $X$ be a Banach space and $K$ a compact space.
\begin{enumerate}
\item If $Y$ is a closed subspace of $X$, then $\ell K(Y)\leq \ell K(X)$.
\item If $Y$ is a total subset of $X$, then $\ell K(X)\leq \ell K(Y)$.
\item If $Y$ is a subset of $C(K)$ which separates the points of $K$, then $\ell K(C(K)) = \ell K(C_p(K))\leq \ell K(Y,\tau_p)$.
\end{enumerate}
\end{prop}

\section{The compact generation index}

In this section we state some main properties of the index $CG(X)$ on
Banach spaces which are analogous to well known properties of weakly
compactly generated spaces. The first observation is that for a Banach
space $X$, $CG(X)$ equals the least infinite number of weakly compact
subsets of $X$ whose union is linearly dense in $X$, $X =
\overline{span}\bigcup_{\lambda<\kappa}K_\lambda$, since if
$\bigcup_{i<\kappa}K_i$ is family of compacta whose union is linearly
dense, we define $K_F = |F|\overline{co}\bigcup_{i\in F}K_i$ for each
finite subset $F$ of $\kappa$, and then we have another family of the same
cardinality whose union is dense. In the following Theorem
\ref{kappaEberlein} we introduce and expose the main properties of
$\kappa$-Eberlein compacta, which generalize well known facts about
Eberlein ($\omega$-Eberlein under this notation) compacta:

\begin{thm}\label{kappaEberlein}
Let $K$ be a compact space and $\kappa$ an infinite cardinal. The
following are equivalent:

\begin{enumerate}

\item $CG(C(K)) \leq \kappa$.

\item $CG(C_p(K)) \leq \kappa$.

\item $K$ is a subspace of a product
$\prod_{\alpha<\kappa}K_\alpha$ of $\kappa$ many factors in which each
$K_\alpha$ is an Eberlein compact.

\item The space can be found as a subset $K\subset [0,1]^\Gamma$ so that
$\Gamma = \bigcup_{\lambda<\kappa}\Gamma_\lambda$ and for every $x\in K$,
and every $\lambda<\kappa$ the set $\{\gamma\in\Gamma_\lambda :
x_\gamma\neq 0\}$ is finite.
\end{enumerate}

When these conditions are satisfied, we say that $K$ is a
$\kappa$-Eberlein compact.
\end{thm}

Proof:  That (1) implies (2) is clear since the weak topology is finer
than the pointwise topology. That (2) implies (3): if
$CG(C_p(K))\leq\kappa$ then we have a family $\{S_\lambda : \lambda<\kappa
\}$ of compact subsets of $C_p(K)$ whose union is dense in $C_p(K)$. In
this situation each compact $K_\lambda$, defined as the quotient of $K$ by
the relation $x\sim y$ iff $f(x)=f(y)$ for all $f\in S_\lambda$, is an
Eberlein compact since $S_\lambda$ is a pointwise compact subset of
continuous functions separates the points of $K_\lambda$ and on the other
hand, $K$ is a subspace of $\prod_{\lambda<\kappa}K_\lambda$. That (3)
implies (4) follows immediately from the well known fact, consequence of
the Amir-Lindenstrauss Theorem \cite{AmiLin}, that $K$ is Eberlein compact
if and only if it verifies (4) for $\kappa=\omega$. For (3) implies (1),
$K$ is a subspace of $ L = \prod_{\alpha<\kappa}K_\alpha$ where each
$K_\alpha$ is Eberlein compact. Since $(C(K),w)$ is a continuous image of
$(C(L),w)$ it is enough to see that $CG(C(L))\leq \kappa$. For each finite
subset $F$ of $\kappa$ we consider $K_F = \prod_{\alpha\in F}K_\alpha$
which is an Eberlein compact. The natural projection $L\To K_F$ induces a
one-to-one operator $T:C(K_F)\To C(L)$, and since $C(K_F)$ is weakly
compactly generated, $CG(T(C(K_F)))=\omega$. The Stone-Weierstrass theorem
implies that $\bigcup\{T(C(K_F)) : F\in [\kappa]^{<\omega}\}$ is dense in
$C(L)$, so $CG(C(L))\leq\kappa$.$\qed$

\begin{thm}\label{boladualkEberlein}
A Banach space $X$ is a subspace of a Banach space $Y$ with
$CG(Y)\leq\kappa$ if and only if $(B_{X^\ast},w^\ast)$ is
$\kappa$-Eberlein.
\end{thm}

PROOF: If $K=(B_{X^\ast},w^\ast)$ is $\kappa$-Eberlein, then
$CG(C(K))\leq\kappa$ and $X$ is a subspace of $C(K)$. Conversely, let
$\{Y_\alpha\}_{\alpha<\kappa}$ be a union of $\kappa$ weakly compactly
generated subspaces of $Y$ which is linearly dense in $Y$. By the Theorem
of Amir and Lindenstrauss \cite{AmiLin}, for every $\alpha$ there is a
one-to-one weak$^\ast$-to-pointwise continuous operator of norm 1,
$T_\alpha:Y_\alpha^\ast\To c_0(\Gamma_\alpha)$ which induces by
composition a weak$^\ast$-to-pointwise continuous operator
$T_\alpha':Y^\ast\To c_0(\Gamma_\alpha)$. Finally, we have a one-to-one
weak$^\ast$-to-pointwise continuous function $\prod T_\alpha: Y^\ast\To
\prod c_0(\Gamma_\alpha)$ which implies that $B_{Y^\ast}$ is a
$\kappa$-Eberlein compact in the weak$^\ast$ topology. Since $B_{X^\ast}$
is a continuous image of $B_{Y^\ast}$ it remains to show that:

\begin{thm}\label{imagencontinuakEberlein}
Any continuous image of a $\kappa$-Eberlein compact is $\kappa$-Eberlein
compact.
\end{thm}

Unlike the previous results, now
Theorem~\ref{imagencontinuakEberlein} cannot be easily deduced
from the well known (and difficult) particular case in which
$\kappa=\omega$. We must instead adapt the arguments of the
countable case to this more general context. We know at least two
ways to do this but since we do not need here this result, and in
any case its proof does not require really new ideas but just
substituting some appearance of the countable power by an
arbitrary $\kappa$, we just indicate how this can be done. One
possibility is to follow the proof of \cite{BenRudWag} changing in
the appropriate places the role of convergent sequences by nets
indicated in the lattice of finite subsets of $\kappa$. For the
other argument, we recall that a $\kappa$-Corson compact is a
compact space which can be found as a subspace $K\subset
\mathbb{R}^\Gamma$ such that for every $x\in  K$,
$|\{\gamma\in\Gamma : x_\gamma\neq 0\}|\leq\kappa$. It is a
consequence of a result of Bell and Marciszewski
\cite{BelMartightness}, who generalized an argument by
Pol~\cite{PolCpCorson}, that the continuous image of a
$\kappa$-Corson compact is $\kappa$-Corson. On the other hand, if
in the terminology of \cite{Namiokanote} we define a
$\kappa$-quasi-RN compact to be a compact whose diagonal is the
intersection of $\kappa$-many almost neighborhoods of the
diagonal, then following analogous arguments as in
\cite{Arvanitakis} and \cite{Namiokanote} it is possible to prove
that the continuous image of a $\kappa$-quasi-RN compact is again
$\kappa$-quasi-RN and that a compact space is $\kappa$-Eberlein if
and only if it is
$\kappa$-Corson and $\kappa$-quasi-RN.\\

We point out that the mentioned result of \cite{BelMartightness} also says
that $\ell(C_p(K))\leq\kappa$ whenever $K$ is a $\kappa$-Corson compact
and, as a consequence, if the dual unit ball of a Banach space $X$ is a
$\kappa$-Corson compact in its weak$^\ast$ topology, then
$\ell(X)\leq\kappa$.

\section{Adequate families}

The examples which we shall present will be based on adequate families of
sets, a concept introduced by Talagrand~\cite{WKA}, precisely to provide
this kind of examples in the countable case. We state in this section, for
the reader's convenience, the facts about this construction that we shall
need. A family of subsets $\mathcal{A}$ of a given set $\Delta$ is called
an adequate family if for every $A$, $A$ belongs to $\mathcal{A}$ if and
only if every finite subset of $A$ belongs to $\mathcal{A}$. Associated to
such a family, we have the compact space
$K_\mathcal{A}\subset\{0,1\}^\Delta$ of the characteristic functions of
elements of $\mathcal{A}$. The key fact proved by Talagrand is:

\begin{thm}\label{uscoadecuados}
If $\Delta$ is a topological space and $\mathcal{A}$ is an adequate family
of subsets of $\Delta$ which consists only of closed subsets of $\Delta$,
then there is an usco map $\phi:\Delta\To 2^X$ with $X$ a point-separating
subset of $C_p(K_\mathcal{A})$.
\end{thm}

Namely $X = \Delta\cup\{0\}$ is the point-separating subset of
$C_p(K_\mathcal{A})$ whose pointwise topology is the narrowest such that
each set $A\in\mathcal{A}$ is closed and the usco map is $\Delta\To 2^X$
given by $\delta\mapsto \{0,\delta\}$. Making use of Proposition
\ref{propiedadesMaria} we get:

\begin{cor}\label{adecuadosindices}
If $\Sigma$ is a complete metric space of weight at most $\kappa$ and
$\mathcal{A}$ is an adequate family of closed subsets of $\Sigma$, then
$\ell K(C(K_\mathcal{A}))\leq\kappa$.
\end{cor}

We will also make use of the following fact:

\begin{thm}\label{weakFarmaki}
Let $\mathcal{A}$ an adequate family of subsets of a set $\Delta$ and
suppose that $K_\mathcal{A}$ is a $\kappa$-Eberlein compact. Then, there
exists a decomposition $\Delta=\bigcup_{\lambda<\kappa}\Delta_\lambda$
such that for every $x\in K_\mathcal{A}$ and every $\lambda<\kappa$, $x$
has only finitely many nonzero coordinates in $\Delta_\lambda$.
\end{thm}

The proof is analogous to that of \cite[Theorem 3.4.2]{FabianWA}, just
changing countable families by families of cardinality at most $\kappa$
where necessary.$\qed$

\section{Weakly Lindelöf determined spaces of arbitrarily high compact generation
index}\label{largeCG}

In the following lemma, we state a standard fact of cardinal
arithmetics which is important for all the discussions afterwards,
namely that there are arbitrarily large cardinals $\tau$ with
$\tau<\tau^\omega$.

\begin{lem}\label{largetau}
Let $\{\kappa_n\}_{n<\omega}$ be a sequence of cardinals such that
$\kappa_{n+1}>\kappa_n$ for every $n<\omega$, and let $\tau$ be
the supremum of this sequence. Then $\tau^\omega>\tau$. In
particular, for any cardinal $\kappa$ there exists a cardinal
$\tau>\kappa$ such that $\tau<\tau^\omega$.
\end{lem}

Proof: Suppose that $\tau^\omega=\tau$ and let $f:\tau\To
\tau^\omega$ be a surjection from $\tau$ onto the set of sequences
of elements of $\tau$, $\alpha\mapsto f(\alpha) =
(f(\alpha)_n)_{n<\omega}$. For every $n<\omega$ the set
$A_n=\{f(\alpha)_n : \alpha<\kappa_n\}$ has cardinality less than
or equal to $\kappa_n<\tau$ hence we can choose
$\beta_n\in\tau\setminus A_n$. The sequence
$\beta=(\beta_n)_{n<\omega}$ is not in the image of $f$ which is a
contradiction. The reason is that if $\beta=f(\alpha)$ for some
$\alpha<\tau$, then there exists $n<\omega$ with $\alpha<\kappa_n$
and therefore $\beta_n = f(\alpha)_n\in A_n$ which is absurd.$\qed$\\

Proof of Theorem~\ref{bigCG}: It is a consequence of Lemma
\ref{largetau} that there exists a metric space $Z$ which cannot
be expressed as the union of $\kappa$ many discrete subsets.
Indeed, if $\kappa<\tau<\tau^\omega$, then $Z=\tau^\omega$ has
this property because the weight of $Z$ is $\tau$ so its discrete
subsets have cardinality less than or equal to $\tau$. From such a
metric space $Z$, we shall construct a Corson compact space $K$
which is not $\kappa$-Eberlein. The construction is in a similar
spirit as those carried out in \cite{ArgMerWLD} to give several
examples of Corson compacta with specific properties. We consider
a well order $<$ on $Z$. Set $\mathcal{A}$ the family of all
subsets $A$ of $Z$ such that each finite subset is of the form
$\{\xi_1<\ldots<\xi_n\}$ with $d(\xi_i,\xi_j)\leq\frac{1}{i}$ for
$i<j$. Notice that this is an adequate family and that every set
$A\in\mathcal{A}$ is either finite or countable. Indeed the order
type of any $A\in\mathcal{A}$ in the well order $<$ cannot be
greater than $\omega+1$ because if
$\{\xi_1<\xi_2<\ldots<\xi_\omega<\xi_{\omega+1}\}$ belongs to
$\mathcal{A}$, then $\xi_\omega = \lim_{n\To\infty}\xi_n =
\xi_{\omega+1}$ which is a contradiction. Hence,
$K=K_\mathcal{A}\subset\{0,1\}^Z$ is a Corson compact and we shall
see that it is not $\kappa$-Eberlein compact. Namely, if it were
$\kappa$-Eberlein, by Theorem \ref{weakFarmaki} there should be a
decomposition $Z=\bigcup_{\lambda<\kappa}Z_\lambda$ such that each
$A\in\mathcal{A}$ has only finitely many elements in each
$Z_\lambda$. We can choose $\lambda$ with $Z_\lambda$ not discrete
and take $z$ an accumulation point of $Z_\lambda$. We shall find
an infinite subset of $Z_\lambda$ which belongs to $\mathcal{A}$
thus getting a contradiction. We take $\xi_1$ the first element of
$Z_\lambda$ such that $\xi_1\in B(z,1)$ (we are denoting by
$B(x,\varepsilon)$ the open ball with center $x$ and radius
$\varepsilon$ in the space $Z$). Second, we take $\xi_2$ the first
element in $Z_\lambda$ greater than $\xi_1$ such that $\xi_2\in
B(z,\frac{1}{2})\cap B(\xi_1,1)$. In the $n$-th step, if
$\xi_1<\cdots<\xi_{n-1}$ have been defined we choose $\xi_n$ to be
the first element of $Z_\lambda$ greater than $\xi_{n-1}$ such
that $\xi_n\in B(z,\frac{1}{n})\cap
\bigcap_{i=1}^{n-1}B(\xi_i,\frac{1}{i})$. After this construction,
$\{\xi_n : n<\omega\}$ is an infinite element
of $\mathcal{A}$ inside $Z_\lambda$.\\

We set $X=C(K)$ with $K$ the compact space defined above. Since $K$ is not $\kappa$-Eberlein, $CG(X)>\kappa$. On
the other hand, it has been proved in \cite[Proposition 4.10]{ArgMerNeg} that if $K$ is a compact subset of
$\{0,1\}^{\alpha}$ such that the order type of the supports of all elements of $K$ is uniformly bounded by a
countable ordinal, then $C(K)$
is weakly Lindel\"of determined.$\qed$\\

Since always $CG(X)\leq 2^{Nag(X)}$ and $CG(X)\leq\ell\Sigma(X)^\omega$,
Theorem \ref{bigCG} also shows that there are weakly Lindelöf determined
Banach spaces of arbitrarily large indices $Nag(X)$ and $\ell\Sigma(X)$.
The Banach spaces such that $\ell\Sigma(X)=\omega$ (or equivalently
$Nag(X)=\omega$) are called weakly countably determined \cite{WKA}. Hence,
for $\kappa\geq 2^\omega$, Theorem~\ref{bigCG} provides examples of weakly
Lindelöf determined Banach spaces which are not weakly countably
determined and of Corson compact spaces which are not Gul'ko compact.

\section{The relation between the compact generation index and the
$\mathcal{K}$-analyticity index}\label{sectrelacionCGLK}

Notice first that if a topological space $Y$ is union of $\tau$
many compacta $\{K_\lambda\}_{\lambda<\tau}$, then $\ell
K(Y)\leq\tau$ because we can get an usco $\phi:\tau\To 2^Y$ by
$\phi(\lambda) = K_\lambda$. Using Proposition
\ref{propiedadesMaria} we get as a consequence that for any Banach
space $X$:

\begin{equation}\label{desigualdadlKCG}
\ell K(X)\leq CG(X)
\end{equation}

On the other hand, if $\ell K(X)\leq\kappa$ then there is an usco
$\phi:\kappa^\omega\To 2^X$ and since the continuous image of a
compact space by an usco is compact and $\kappa^\omega$ is the
union of $\mathbf{d}(\kappa)$ many compacta,

\begin{equation}\label{desigualdadCGlK}
CG(X)\leq \mathbf{d}(\ell K(X))
\end{equation}

Finally, the last relation is that for any Banach space $X$,

\begin{equation}\label{desigualdaddenslKCG}
dens(X)\geq\mathbf{b}_{\ell K(X)}(CG(X))
\end{equation}

We prove it by contradiction. Suppose the contrary and call
$\delta=dens(X)$, $\kappa=\ell K(X)$ and $\tau =CG(X)$ so that
$\delta<\mathbf{b}_\kappa(\tau)$. We have an usco $\phi:\kappa^\omega\To
2^X$ and we can find a subset $\Sigma\subset\kappa^\omega$ of cardinality
$\delta$ such that $\phi(\Sigma)$ is dense in $X$. Since
$\delta<\mathbf{b}_\kappa(\tau)$, $\Sigma$ is a subset of a union of less
than $\tau$ many compact subsets of $\kappa^\omega$, so $CG(X)<\tau$, a
contradiction.\\

Relations (\ref{desigualdadlKCG}) - (\ref{desigualdaddenslKCG}) already
prove one implication of Theorem \ref{relacionCGWK}. Before passing to the
converse, we make an observation about the evaluation of the indices on a
generalized Cantor cube: For $K=\{0,1\}^\kappa$ and $X=C(K)$ we have that
$\ell(X) = \ell K(X)= CG(X) = \kappa$. On the one hand, clearly
$\{0,1\}^\kappa$ is $\kappa$-Eberlein. On the other hand, the evaluation
maps $D = \{\delta_x : x\in\kappa\}$ constitute a discrete pointwise
closed subset of $C(\{0,1\}^\kappa)$, so $\ell(D)=\kappa$ and
$\ell(X)\geq\kappa$.\\

Now we fix cardinals $\kappa$, $\tau$ and $\delta$ like in part (1) of
Theorem \ref{relacionCGWK} and we will show a Banach space like in part
(2). First, we take $S$ a subset of $\kappa^\omega$ of cardinality
$\mathbf{b}_\kappa(\tau)$ which can be decomposed into $\tau$ many pieces
$S = \bigcup_{\lambda<\tau}S_\lambda$ verifying the following two
properties:
\begin{itemize}
\item[(S.1)] $S_\lambda$ is not contained in any union of less than $\tau$
many compacta of $\kappa^\omega$ \item[(S.2)] There is a subset $U\subset
S$ of cardinality $\kappa$ such that $|U\cap S_\lambda|\leq 1$ for all
$\lambda$ and such that $x_0\neq y_0$ for any two different elements $x,
y\in U$ .
\end{itemize}

We can construct such an $S$ as follows: Take $A$ a subset of
$\kappa^\omega$ of cardinality $\mathbf{b}_\kappa(\tau)$ which cannot be
covered by less than $\tau$ many compacta of $\kappa^\omega$ and
$A'=\{a_\lambda : \lambda<\tau\}$ a subset of $A$ of cardinality $\tau$.
Without loss of generality, we suppose that the set $U
=\{x\in\kappa^\omega : x_0=x_n \forall n<\omega\}$ of the constant
sequences is a subset of $A'$. For $\lambda<\tau$, we define $S_\lambda =
\{x\in\kappa^\omega : x_{2n} = (a_\lambda)_n,\
(x_{2n+1})_{n\in\omega}\in A\}$ and $S=\bigcup_{\lambda<\tau}S_\lambda$.\\

We construct now compact space inspired on the example of Talagrand
\cite{WKA} of a Talagrand non Eberlein compact. Consider $\mathcal{A}$ the
family of all subsets $A\subset S$ which verify the two following
properties:
\begin{itemize}
\item[(A.1)] There exists some $n(A)<\omega$ such that for every different
elements $x,y\in A$, it happens that $x_{n(A)}\neq y_{n(A)}$ but $x_m=y_m$
for all $m<n(A)$.
\item[(A.2)] $|A\cap S_\lambda|\leq 1$ for every $\lambda<\tau$.
\end{itemize}

This $\mathcal{A}$ is an adequate family which consists of closed subsets
of $S$, hence the compact space $L = K_\mathcal{A}\subset \{0,1\}^S$ is a
$\kappa$-Talagrand compact, by Corollary~\ref{adecuadosindices}. We define
$Y=C(L)$, so that $\ell K(Y)\leq\kappa$. Indeed, $\ell K(Y)=\kappa$
because if $U$ is a set like in (S.2) then all subsets of $U$ belong to
$\mathcal{A}$ and there is therefore a copy of $\{0,1\}^\kappa$ inside $L$
and, as we observed, $\ell K(C(\{0,1\}^\kappa)=\kappa$.\\

On the other hand, property (A.2) implies that $L\subset \{0,1\}^S$ is a
$\tau$-Eberlein compact (indeed the partition $S =
\bigcup_{\lambda<\tau}S_\lambda$ fulfills the conditions in Theorem
\ref{kappaEberlein}(4)) and hence, $CG(Y)\leq\tau$. We check now that
precisely $CG(Y)=\tau$. Assume by contradiction that $L$ is
$\tau'$-Eberlein for some $\tau'<\tau$. Then by Lemma~\ref{weakFarmaki} we
would find a partition $S = \bigcup_{i<\tau'}\Delta_i$ such that each
element of $L$ has only finitely many nonzero coordinates in each
$\Delta_i$. Let us analyze for a moment what this condition on $\Delta_i$
means. For a subset $F$ of $\kappa^n$, we denote $F\times\kappa^{>n} =
\{x\in\kappa^\omega : (x_0,\ldots,x_{n-1})\in F\}$ and if $G\subset
\kappa^m$ with $m>n$ we write $F<G$ if the restrictions of the elements of
$G$ to the first $n$ coordinates constitute precisely the set $F$. The
fact that we cannot find an infinite subset $A$ of $\Delta_i$ satisfying
(A.2) and also (A.1) with n(A)=0 implies that there exist finite sets
$F_0\subset\kappa$ and $G_0\subset\tau$ such that $\Delta_i\subset
F_0\times\kappa^{>0}\cup\bigcup_{\lambda\in G_0}S_\lambda$. Analogously,
paying attention in each step to sets $A$ with $n(A)=n$, we can find
inductively for every $n<\omega$ finite sets $F_n\subset\kappa^n$ and
$G_n\subset \tau$ such that $F_{n-1}<F_n$, $G_{n-1}\subset G_n$ and
$\Delta_i\subset F_n\times\kappa^{>n}\cup\bigcup_{\lambda\in
G_n}S_\lambda$. This implies that for every $i<\tau'$ we can find a
compact set $K_i = \bigcap_{n<\omega}F_n\times\kappa^{>n}$ of
$\kappa^\omega$ and a countable set $G_i\subset\tau$ such that $\Delta_i
\subset K_i\cup\bigcup_{\lambda\in G_i}S_\lambda$. Hence,
$$S = \bigcup_{i<\tau'}K_i\cup\bigcup_{\lambda\in\bigcup_{i<\tau'}
G_i}S_\lambda.$$ Since $|\bigcup_{i<\tau'} G_i|\leq\tau'\cdot\omega<\tau$,
we can take $\lambda_0\not\in\bigcup_{i<\tau'}G_i$ and then
$S_{\lambda_0}$ is covered by $\tau'<\tau$ compact subsets of
$\kappa^\omega$, a contradiction.\\

So far, we know that $\ell K(Y) = \kappa$ and $CG(Y)=\tau$. Since
$|S|=\mathbf{b}_\kappa(\tau)$, this is the weight of $L\subset\{0,1\}^S$
(it could not be lower because of the general relation
(\ref{desigualdaddenslKCG})), hence $dens(Y) = \mathbf{b}_\kappa(\tau)$.\\

Finally, we consider the space $X=C(L)\oplus c_0(\delta)$. The space
$c_0(\delta)$ is weakly compactly generated, so
$CG(c_0(\delta))=\omega=\ell K(c_0(\delta))$ and
$dens(c_0(\delta))=\delta$. All these indices, when considered on a finite
product, take as value the maximum of the value of each factor, so $X$ is
the space we were looking for.

\section{The number of compact spaces which generate a
subspace}\label{sectsubespacios}

This section is devoted to the proof of Theorem~\ref{thCGsubespacios}. The
situation is very similar to the previous section. Let us assume that we
are in the situation of part (2) of Theorem \ref{thCGsubespacios}. That
$\tau\leq\delta$ is evident. Being $Y$ a closed subspace of $X$ then $\ell
K(Y)\leq \ell K(X)$, and from this and (\ref{desigualdadCGlK}), we have
that
$$\tau = CG(Y) \leq\mathbf{d}(\ell K(Y))\leq\mathbf{d}(\ell K(X))
=\mathbf{d}(\kappa),$$ and since $\ell K(Y)\leq\kappa$, by
(\ref{desigualdaddenslKCG}),

$$\delta = dens(Y)\geq \mathbf{b}_{\ell K(Y)}(CG(Y)) = \mathbf{b}_{\ell K(Y)}(\tau)\geq
\mathbf{b}_{\kappa}(\tau).$$\

For the converse, if $\tau<\min(\kappa,\delta)$, it is enough to take
$X=C(\{0,1\}^\kappa)\oplus c_0(\delta)$ and $Y=C(\{0,1\}^\tau)\oplus
c_0(\delta)$. So we assume from now on that
$\kappa\leq\tau\leq\mathbf{d}(\kappa)$ and
$\delta\geq\mathbf{b}_\kappa(\tau)$, and we will adapt an example of
Argyros \cite[Section 1.6]{FabianWA} by similar modifications as in the
proof of Theorem \ref{relacionCGWK}. First, we take like in that proof a
subset $S$ of $\kappa^\omega$ of cardinality $\mathbf{b}_\kappa(\tau)$
which can be decomposed into $\tau$ many pieces $S =
\bigcup_{\lambda<\tau}S_\lambda$ verifying (S.1) and (S.2). We consider
the compact space $K\subset [0,1]^{S}$ which consists of the function of
the form $\frac{1}{n}\chi_A$ ($\chi_A$ is the characteristic function of
the set $A$) for some natural number $n$ and some set $A$ satisfying:

\begin{itemize}
\item[(B.1)] For every different elements $x,y\in A$, it happens
that $x_{n}\neq y_{n}$ but $x_m=y_m$ for all $m<n$.
\item[(B.2)] $|A\cap S_\lambda|\leq 1$ for every $\lambda<\tau$.
\end{itemize}

We note that the decomposition $S= \bigcup_{t\in\kappa^n}\{\sigma\in S :
\sigma|_n = t\}$ verifies the conditions of Theorem \ref{kappaEberlein}(4)
for $\varepsilon=\frac{1}{n}$, so $K$ is a $\kappa$-Eberlein compact and
since, again, it contains a copy of $\{0,1\}^\kappa$, $CG(C(K))=\kappa$.
For every $\sigma\in S$ we consider the ``projection'' function
$f_\sigma\in C(K)$ and we set $Y$ the subspace generated by $\{f_\sigma :
\sigma\in S\}$. For every $\lambda<\tau$, it follows from condition (B.2),
that $\{f_\sigma : \sigma\in S_\lambda\}\cup\{0\}$ is a pointwise (hence
weakly) compact subset of $C(K)$. Therefore, $CG(Y)\leq\tau$. We suppose
by contradiction that $CG(Y)=\theta<\tau$. The rest of the proof follows
closely that of \cite[Theorem 1.6.3]{FabianWA}. As a consequence of the
Amir-Lindenstrauss Theorem, we can find then a set of generators of $Y$
like $\{y_\delta : \delta\in\Delta\}$ such that $\Delta =
\bigcup_{\eta<\theta}\Delta^\eta$ and $\{y_\delta :
\delta\in\Delta^\eta\}\cup\{0\}$ is homeomorphic in the weak topology to
the one compactification of a discrete set, being 0 the point of infinity.
We define a function $F:S\times\Delta\To\mathbb{R}$ like $F(\sigma,\delta)
= y_\delta(\chi_{\{\sigma\}})$.\\

Statement 1: For every $\sigma\in S$, $1\leq |\{\delta\in\Delta :
F(\sigma,\delta)\neq 0\}|\leq\theta$. Namely, if that set were empty,
since $\{y_\delta : \delta\in\Delta\}$ generates $Y$ that would mean that
$y(\chi_{\{\sigma\}})=0$ for all $y\in Y$ which is false for $y=f_\sigma$.
On the other hand, since 0 is the ``weak limit point'' of $\{y_\delta :
\delta\in\Delta^\eta\}$ each set $\{\delta\in\Delta^\eta :
|F(\sigma,\delta)|>\frac{1}{m}\}$ is finite.\\

Statement 2: For every $\delta\in\Delta$, $|\{\sigma\in S :
F(\sigma,\delta)\neq 0\}|\leq\omega$. Indeed, each $\{\sigma\in S :
|F(\sigma,\delta)|>\frac{1}{m}\}$ is finite because we can find an element
$y\in Y$ which is a linear combination of some
$f_{\sigma^1},\ldots,f_{\sigma^k}$ with $\|y-y_\delta\|<\frac{1}{2m}$ and
in this case, whenever $|F(\sigma,\delta)|>\frac{1}{m}$, it must be the
case that $\sigma=\sigma^i$ for some $i\leq k$.\\

From these two statements, playing back and forth we find a partition $S
=\bigcup_{\alpha<\lambda}\Gamma_\alpha$ and disjoint sets $\Delta_\alpha$,
all sets of cardinality at most $\theta$ such that whenever
$F(\sigma,\delta)\neq 0$ then there exists $\alpha<\lambda$ such that
$\sigma\in\Gamma_\alpha$ and $\delta\in\Delta_\alpha$.\\

Since $|\Gamma_\alpha|\leq\theta$ we can enumerate it as
$\{\sigma_\nu^\alpha : \nu<\theta\}$ (with repetitions if necessary).
Then, we set $\Sigma_\nu = \{\sigma^\alpha_\nu : \alpha<\lambda\}$ so that
$S = \bigcup_{\nu<\theta}\Sigma_\nu$. Further, for $\nu<\theta$,
$m\in\mathbb{N}$ and $\eta<\theta$ we set
$$\Sigma_{\nu m\eta} = \left\{\sigma\in\Gamma_\nu : \exists
\delta\in\Delta^\eta : y_\delta(\chi_{\{\sigma\}})>\frac{1}{m}\right\}.$$
By statement 1, $\Sigma_\nu = \bigcup_{m\eta}\Gamma_{\nu m\eta}$ and
moreover $S = \bigcup_{\nu<\theta,m<\omega,\eta<\theta}\Sigma_{\nu
m\eta}$.\\

We proved in the previous section that for such a decomposition of
$S$ into $\theta<\tau$ pieces, there must exist $\nu$, $m$ and
$\eta$ such that $\Sigma_{\nu m\eta}$ is not contained in a
compact subset of $\kappa^\omega$, and this implies that, for some
$n<\omega$ there is an infinite set $A$ satisfying (B.1) and (B.2)
which is contained in $\Sigma_{\nu m\eta}$, call it $A =
\{\sigma_i : i<\omega\}$. Call $\alpha_i$ the only ordinal such
that $\sigma_i\in\Gamma_{\alpha_i}$. Since $A\subset\Sigma_\nu$,
$\sigma_i = \sigma_\nu^{\alpha_i}$ and if $i\neq j$ then
$\alpha_i\neq\alpha_j$. Also, since $A\subset \Sigma_{\nu m\eta}$
for every $i$, there exists $\delta_i\in\Delta^\eta$ such that
$y_{\delta_i}(\chi_{\{\sigma_i\}})>\frac{1}{m}$. Notice that if
$i\neq j$ then, $\delta_i\in\Delta_{\alpha_i}$ so
$\delta_i\not\in\Delta_{\alpha_j}$
and $y_{\delta_i}(\chi_{\{\sigma_j\}})=0$.\\

Let now $B$ be a finite subset of $A$. Then, for every
$\sigma\in\Sigma$ we have
$$f_\sigma(\frac{1}{l}\chi_B) = \frac{1}{l}\chi_B(\sigma) =
\frac{1}{l}\sum_{\sigma'\in B}\chi_{\{\sigma'\}}(\sigma) =
\frac{1}{l}\sum_{\sigma'\in B}f_\sigma(\chi_{\{\sigma'\}}).$$

Hence, $y(\frac{1}{l}\chi_B) = \frac{1}{l}\sum_{\sigma'\in
B}y(\chi_{\{\sigma'\}})$ for every $y\in Y$. Now let $B_1\subset
B_2\subset\cdots$ be a sequence of finite subsets of $A$ whose
union is $A$. Then $\frac{1}{l}\chi_{B_j}\To\frac{1}{l}\chi_A$ and
so

$$|y_{\delta_i}(\frac{1}{l}\chi_A)| =
\lim_{j\To\infty}|y_{\delta_i}(\frac{1}{l}\chi_{B_j})| =
\frac{1}{l}|y_{\delta_i}(\chi_{\{\sigma_i\}})|>\frac{1}{lm}$$

This is a contradiction since $\{\delta_i : i<\omega\}\subset\Delta^\eta$
so it weakly converges to 0.$\qed$

\section{Remarks on other indices}\label{otherindices}

The only example that we know of a Banach space $X$ in which
$\ell\Sigma(X)<\ell K(X)$ is one of Talagrand \cite{WCD}, a variation of
which can be taken of density character $\omega_1$, cf. \cite{cardinalb}.
In this case $\ell\Sigma(X)=\omega<\omega_1 = \ell K(X)$. The
``enlargement'' of such an example offers a number of difficulties and we
do not even know whether there exists some Banach space $X$ with
$\omega<\ell\Sigma(X)<\ell K(X)$.\\

About the Nagami index, it is provided in \cite{CasMunOri} an example of a
topological space $Y$ with $Nag(Y)<\ell\Sigma(Y)$. We shall provide next
examples of this kind for spaces of the form $Y=C_p(K)$ with $K$ compact.
This does not provide yet examples of Banach spaces for which the two
indices do not coincide, because it is not clear whether the Nagami index
coincides for the weak and the pointwise topology: the proof of the fact
that $\ell\Sigma (C(K)) = \ell\Sigma(C_p(K))$ \cite{CasMunOri} heavily
depends on the
fact that closure points in metric spaces are limits of sequences.\\

\begin{thm}\label{NagamiCG}
Let $\kappa$ be any infinite cardinal. Then there exists a compact space
$K$ with $Nag(C_p(K))\leq\kappa<CG(C_p(K))$.
\end{thm}

We already observed that we always have $CG(X)\leq \ell\Sigma(X)^\omega$,
so when $\kappa=\kappa^\omega$ the compactum of Theorem \ref{NagamiCG}
verifies
$Nag(C_p(K))<\ell\Sigma(C_p(K))$.\\

Proof: Consider the space $T=\kappa^\kappa$ (the product of
$\kappa$ many discrete spaces of size $\kappa$) which has weight
$\kappa$. We consider the adequate family $\mathcal{A}$ of all
subsets $A$ of $T$ such that there is $\lambda<\kappa$ such that
for any $x\neq y$ in $A$, $x|_{[0,\lambda)}=y|_{[0,\lambda)}$ and
$x_\lambda\neq y_\lambda$. We take $K=K_{\mathcal{A}}$. We know
from Theorem~\ref{uscoadecuados} that there is an usco $T\To
2^{Y}$ with $Y$ a subset of $C_p(K)$ that separates the points of
$K$. In an analogous way as it is proven for the index
$\ell\Sigma(X)$ \cite{CasMunOri}, this implies that
$Nag(C_p(K))\leq\kappa$. Suppose now that $K$ were
$\kappa$-Eberlein. Then, by Theorem~\ref{weakFarmaki}, we should
be able to find a decomposition $T=\bigcup_{i<\kappa}T_i$ such
that any set of $\mathcal{A}$ has only finitely many elements in
each $T_i$. We will find $t\in T=\kappa^\kappa$ such that
$t\not\in T_i$ for any $i<\kappa$, thus obtaining a contradiction.
We define it inductively. The set $\{x_0 : x\in
T_0\}\subset\kappa$ is finite so we may take $t_0$ out of it. This
will guarantee that $t\not\in T_0$. If we already defined $t_j$
for $j<i$, the set $\{x_i : x\in T_i\text{ and }x_j=t_j \ \forall
j<i\}$ is finite, so we can choose $t_i$ out of it. This
guarantees that $t\not\in T_i$.$\qed$\\

\end{document}